\documentclass[review]{elsarticle}

\usepackage{lineno,hyperref}
\modulolinenumbers[5]

\journal
{arXiv}
\usepackage{amscd}
\usepackage{mathrsfs}
\usepackage{amsfonts}
\usepackage{graphicx}
\usepackage{amsmath,amscd,stmaryrd,amssymb,array, amsfonts,amsmath,amssymb}
\usepackage{enumitem}

\usepackage{caption}
\usepackage{subfigure}

\numberwithin{figure}{section}
 \numberwithin{equation}{section}

\newtheorem{theorem}{Theorem}[section]
\newtheorem{proposition}[theorem]{Proposition}
\newtheorem{definition}[theorem]{Definition}

\newtheorem{lemma}[theorem]{Lemma}
\newtheorem{remark}[theorem]{Remark}

\newcommand{\bi}{{\mathbf{i}}}

\newcommand{\bbA}{{\mathbb A}}

\newcommand{\bbC}{{\mathbb C}}
\newcommand{\bbH}{{\mathbb H}}

\newcommand{\bbS}{{\mathbb S}}

\newcommand{\cC}{{\mathcal C}}
\newcommand{\cD}{{\mathcal D}}

\newcommand{\cW}{{\mathcal W}}

\newcommand{\sL}{{\mathscr L}}

\def\be{\begin{equation}}
\def\ee{\end{equation}}
\def\bes{\begin{equation*}}
\def\ees{\end{equation*}}
\def\bsp{\begin{split}}
\def\esp{\end{split}}

\def\ba{\begin{array}}
\def\ea{\end{array}}
\def\benu{\begin{enumerate}}
\def\eenu{\end{enumerate}}
\def\bt{\begin{theorem}}
\def\et{\end{theorem}}
\def\bp{\begin{proposition}}
\def\ep{\end{proposition}}
\def\bl{\begin{lemma}}
\def\el{\end{lemma}}
\def\br{\begin{remark}}
\def\er{\end{remark}}
\def\bd{\begin{definition}}
\def\ed{\end{definition}}

\def\b{\beta}

\def\de{\delta}
\def\d{\partial}
\def\pa{\partial}

\def\lam{\lambda}

\def\o{\stackrel{\circ}}
\def\ve{\varepsilon}
\def\sig{\sigma}

\def\a{\alpha}
\def\W{\Omega}

\def\.{\cdot}

\def\bbA{\mathbb{A}}
\def\bbC{\mathbb{C}}

\def\R{\mathbb{R}}

\def\bbX{\mathbb{X}}

\def\A{\forall}
\def\ol{\overline}

\def\ra{\rightarrow}

\def\~{\tilde}
\def\8{\infty}
\def\X{\times}
\def\({\left(}
\def\){\right)}

\def\mb{\mbox}
\def\emp{\emptyset}
\def\sm{\setminus}

\def\bx{\blacksquare}

\def\Hs{\hspace{1cm}}\def\hs{\hspace{0.5cm}}
\def\Vs{\vskip8pt}\def\vs{\vskip4pt}

\def\({\left(}\def\){\right)}

\begin{document}

\begin{frontmatter}

\title{A Note on the Krein-Rutman Theorem for\\ Sectorial Operators$^\dag$\footnote{$\dag$This work was supported by the National Natural Science Foundation of China [11871368]}}


\author[mymainaddress]{Desheng  Li\,\corref{mycorrespondingauthor}}
\cortext[mycorrespondingauthor]{Corresponding author}
\ead{lidsmath@tju.edu.cn}

\author[mymainaddress]{Ruijing Wang}
\ead{wrj\_math@tju.edu.cn}
\address[mymainaddress]{School of Mathematics,  Tianjin University, Tianjin 300350,  China}

\author[mysecondaryaddress]{Luyan Zhou}
\ead{zhouly@bnu.edu.cn}

\address[mysecondaryaddress]{School of Mathematical Sciences, Laboratory of Mathematics and Complex Systems, MOE, Beijing Normal
University, Beijing 100875, China}

\begin{abstract}
In this note we present some generalized versions of the  Krein-Rutman theorem for   sectorial  operators. They are formulated   in  a fashion  that can be easily applied to elliptic operators. Another feature of these generalized versions is that they contain some information on the generalized eigenspaces associated with non-principal eigenvalues, which are helpful in the study of the dynamics of evolution equations in ordered Banach spaces.
\end{abstract}

\end{frontmatter}


\section{Introduction}

It is well known that the Krein-Rutman theorem (see e.g. \cite{Krein, Schab}) for bounded operators (particularly for compact operators) plays a crucial role in the discussion   of  the principal eigenvalue problem of  elliptic operators via their resolvents.  But generally only  part of the information concerning  the principal eigenvalue and eigenvectors  can be obtained. One reason is  that there may  be no one-one correspondence between the boundary spectrum of an  elliptic operator $A$ and the peripheral spectrum of its resolvent $R_\lam(A)$. As a remedy, one has to do some tricky  PDE argument when refined information on the principal eigenvalue and eigenvectors is needed; see e.g. Du \cite[Theorem 1.4]{Du}, Evans \cite[Sec. 6.5, Theorem 3]{Evans}  and Ni \cite{Ni} etc.
It is therefore of particular interest to give variants of  the classical Krein-Rutman theorem for  unbounded operators.

In Greiner et al. \cite{Grein} it was proved that the spectral bound $s(A)$ of the generator $A$ of a positive $C_0$-semigroup on a Banach space with a normal reproducing cone is contained in $\sig(A)$. The monograph  \cite{Arendt} contains a far reaching theory concerning the above question in the framework of generators of $C_0$-semigroups on the functional space $C_0(X)$ consisting of continuous functions vanishing at infinity, where $X$ is a locally compact space; see \cite[Chap. B-III]{Arendt}. 
These extensions  make the  Krein-Rutman theorem to be more efficient  in studying the spectral properties of unbounded  operators.

 When dealing with an elliptic differential operator $L$ on a domain $\W\subset \R^n$ as in \eqref{equ:1-2}, the frequently used  spaces in which $L$ may  generate a  $C_0$-semigroup are the Sobolev ones. However, in many cases the cone  of nonnegative functions can  neither  have interior nor be normal  in a Sobolev space. (The above two properties were   required in most abstract results in the literature.) Of course, $L$ can be considered as an operator in the  space   $C^{k,\a}(\ol\W)$, in which  the cone $\cC$ of nonnegative functions have nonempty interior. But as it was pointed in Kielhofer \cite{Kiel},
$L$ can fail to generate a $C_0$-semigroup  in this space. To cover this situation, Nussbaum \cite{Nuss2} made an effort to  extend  a major part of the Krein-Rutman theorem to operators which may not  generate a  $C_0$-semigroup.

In a recent paper \cite{LiDS}, the authors reformulated the classical Perron-Frobenius theory and Krein-Rutman theorem by using an elementary dynamical approach. Inspired by this work, in this note  we present another  generalization of the    Krein-Rutman theorem for   sectorial  operators in  a  formalism that seems to be more natural and suitable for elliptic operators. It is  different in some ways from those mentioned above and   can be easily  applied to elliptic operators, which   enables  us to reduce significantly the technical  PDE argument  involved in the investigation of elliptic eigenvalue problems. Another feature of these generalized versions is that they contain some information on the generalized eigenspaces associated with non-principal eigenvalues, which are helpful in the study of the dynamics of evolution equations in ordered Banach spaces (see e.g. \cite{LiDS3}). As an illustrating example,  the principal eigenvalue problem of a general elliptic operator associated with degenerate mixed  boundary condition  is discussed. 

The remaining part of this  paper is organized as follows.  In Section \ref{s:2} we do   some preliminary work, and in Section \ref{s:3} we prove generalized   Krein-Rutman type theorems for   sectorial  operators. Section \ref{s:4} consists of an example  illustrating the theoretical results mentioned above.

\section{Preliminaries}\label{s:2}
Let   $X$ be  a real Banach space with norm $\|\.\|$.
Given  $M\subset X$, the {\em interior} and {\em closure } of $M$  are denoted respectively by  $\o M$ and $\ol M$. When we need to emphasize in which space the interior and closure are taken, we also use the notations $\mb{int}_X M$ and $\mb{Cl}_X M$ in place of  $\o M$ and $\ol M$, respectively.
 For  $x\in X$, set $$d(x,M)=\inf_{y\in M}\|x-y\|.$$

\subsection{Some basic knowledge in the spectral theory  of operators}\label{s:2.3}
\vs Let  $A$ be a closed densely defined operator in $X$. Denote by $\sig(A)$ and $\rho(A)$ the {\em spectrum} and {\em resolvent set} of $A$, respectively.  For $\lam\in\rho(A)$, let $R_\lam(A):=(\lam-A)^{-1}$ be the resolvent of $A$.

Let $\bbX=X+\bi X$ be the complexification  of  $X$ (see e.g. \cite[Sect. 7]{LiDS} for details), and define the complexification $\bbA$ of $A$ as
$$
\bbA u=Ax+\bi Ay,\Hs\A\,u=x+\bi y\in \bbX.
$$
Given  $\mu\in\sig( A)=\sig(\bbA)$, let
\be\label{e:es}
\mb{GE}_\mu( \bbA)=\{\xi\in\bbX:\,\,(\bbA-\mu)^{j}\xi=0\mb{ for some }j\geq 1\}.
\ee
Then $\mb{GE}_\mu( \bbA)$ is an invariant subspace of $\bbA$. For each  $\xi\in \mb{GE}_\mu( \bbA)\sm\{0\}$, it is clear that  there is an  integer $k\geq 1$ such that
\be\label{e:2.5}( \bbA-\mu)^{j}\xi\ne0\,\,\,(0\leq j\leq k-1),\hs ( \bbA-\mu)^{k}\xi=0.\ee
For convenience in statement, we call the number $k$  in \eqref{e:2.5}  the {\em rank} of $\xi$, denoted by  $\mb{rank}(\xi)$.

For $\mu\in\sig(A)$, it follows by the invariance of  $\mb{GE}_\mu( \bbA)$ that
\be\label{e:GE}\ba{ll}
         \mb{GE}_\mu(A):=\{\mb{Re}\,\xi:\,\,\xi\in \mb{GE}_\mu(\bbA)\}\,\,\(=\{\mb{Im}\,\xi:\,\,\xi\in \mb{GE}_\mu(\bbA)\}\)
                \ea
         \ee
is an invariant subspace of $A$ in $X$, which will be referred to as the {\em generalized eigenspace} of $A$ pertaining to $\mu$.  (The second equality in \eqref{e:GE} is due to the fact that if $\xi\in \mb{GE}_\mu(\bbA)$, then $\pm \mathbf{i}\xi\in \mb{GE}_\mu(\bbA)$.)

  \bl\label{l:2.6}  Let $\mu\in\sig(A)$, and  $\lam\in \rho(A)$. Then
$$
R_\lam(A)\mb{\em GE}_\mu( A)= \mb{\em GE}_\mu(A)=\mb{\em GE}_{(\lam-\mu)^{-1}}(R_\lam(A)).
$$
\el
{\bf Proof.} It suffices to check that
\be\label{e:2.6}
R_\lam(\bbA)\mb{GE}_\mu( \bbA)= \mb{GE}_\mu(\bbA)=\mb{GE}_{\lam_\mu}(R_\lam(\bbA)),
\ee
where $\lam_\mu:=(\lam-\mu)^{-1}$.

Let $\xi\in \mb{GE}_\mu(\bbA)$. Then $(\bbA-\mu)^k\xi=0$ for some $k\geq 1$. Hence
$$
(\bbA-\mu)^k\(R_\lam(\bbA)\xi\)=R_\lam(\bbA)(\bbA-\mu)^k\xi=R_\lam(\bbA)0=0.
$$
It follows  that $R_\lam(\bbA)\xi\in \mb{GE}_\mu(\bbA)$. Simple calculations also yield
\be\label{e:2.8}
(\lam-\bbA)^k\(R_\lam(\bbA)-\lam_\mu\)^k\xi=\lam_\mu^k(\bbA-\mu)^k\xi=0.
\ee
 Since $\lam\in\rho(\bbA)$, \eqref{e:2.8}  implies   that $\(R_\lam(\bbA)-\lam_\mu\)^k\xi=0$. Therefore $\xi\in \mb{GE}_{\lam_\mu}(R_\lam(\bbA))$.
In conclusion, we have $$R_\lam(\bbA)\mb{GE}_\mu( \bbA)\subset \mb{GE}_\mu(\bbA)\subset\mb{GE}_{\lam_\mu}(R_\lam(\bbA)).$$

The verification of the inverse inclusions is similar. We omit it.
$\bx$

\Vs
Denote by  $\sig_e(A)$     the {\em essential spectrum}   of $A$ in the terminology of Browder \cite[pp. 107-108, Def. 11]{Browder}. Then  each $\mu\in\sig(A)\sm\sig_e(A)$ is isolated in $\sig(A)$ with $\mb{GE}_\mu(A)$ being a finite-dimensional subspace of $X$; see \cite[pp. 108]{Browder}.

           The {\em spectral bound}  $\mb{spb}(A)$ and  {\em essential spectral bound} $\mb{spb}_e(A)$ of $A$ are defined as
$$
\mb{spb}(A)=\sup\{\mb{Re}\,\mu:\,\,\mu\in\sig(A)\},\hs \mb{spb}_e(A)=\sup\{\mb{Re}\,\mu:\,\,\mu\in\sig_e(A)\}.
$$
(We assign $\mb{spb}_e(A)=-\8$ if $\sig_e(A)=\emp$.) Set
 $$\sig_b(A)=\sig(A)\cap\{\mb{Re}\,z=\mb{spb}(A)\}.$$
 $\sig_b(A)$ is called  the {\em boundary spectrum} of $A$.

Let  $\sL(X)$ be the space of bounded linear  operators on $X$.
If $A\in\sL(X)$, we define   the {\em spectral radius} $r(A)$ and  {\em essential spectral radius} $r_e(A)$ as
$$r(A)=\sup\{|\mu|:\,\,\mu\in\sig(A)\},\hs r_e(A)=\sup\{|\mu|:\,\,\mu\in\sig_e(A)\}.$$
It is a basic knowledge   that $r(A)=\lim_{k\ra\8}\|A^k\|^{1/k}$.

\subsection{Cones and positive operators}
\vs
Let $X$ be a Banach space. A  {\em wedge} in $X$ is a closed subset  $ K\subset X$ with $K\ne\{0\}$ such that  $t K\subset  K$ for all $t\geq 0$.
  A convex wedge $K$ with $ K\cap (- K)=\{0\}$ is called a {\em cone}.

A cone   $ K$ is said to be   {\em total} (resp., {\em solid}),  if $\ol{K-K}=X$ (resp., $\o K\ne\emp$).
\vs

From now on we always assume that there has been given a cone $K$ in $X$.


An operator  $A\in \sL(A)$  is called {\em positive} (resp. {\em strongly positive}) if $AK\subset K$ (resp. $A(K\sm\{0\})\subset \o K$).

 \bd \cite{LiDS} In case  $K$ is  a solid cone, a positive operator  $A$ is called  weakly irreducible, if the boundary $\pa K$ of $K$   contains no   eigenvectors of $A$ pertaining to nonnegative  real eigenvalues.
\ed

\br It is almost obvious  that strongly positive operators and primitive operators are weakly irreducible.
 (Recall that $A$ is said to be {\em primitive}, if there is an integer $m\geq 1$ such that $A^m(K\sm\{0\})\subset \o K$; see \cite[pp. 285]{LN}.)
 We also infer from the argument following Definition 7.5 in \cite{LiDS}  that  irreducibility  (in the terminology of \cite{LN})  implies  weak irreducibility defined as above.\er

We finally recall some  generalized  versions of the  Krein-Rutman  Theorem  for bounded operators given  in \cite{LiDS} to conclude this section.

\bt\label{t:7.4}{\cite[Theorem 7.3]{LiDS}}\,   Let $A\in\sL(X)$ be a positive operator with
 \be\label{e:2.6}r_e:=r_e(A)<r(A):=r.\ee
   If $K$ is total, then the following assertions hold true:
        \benu
        \item[$(1)$] $r$ is an eigenvalue of $A$ with  a principal  eigenvector  $u\in K$.
\item[$(2)$] If $\o K$  contains a principal eigenvector $v$ of $A$, then the algebraic and  geometric multiplicities  of $r$ coincide.
              \item[$(3)$] If $\mu$ is a complex eigenvalue with $|\mu|>r_e$, then $\mb{\em GE}_\mu(A)\cap K=\{0\}.$
                  \item[$(4)$] All eigenvectors of $A$ pertaining to eigenvalues  $\mu\ne r$ with $|\mu|>r_e$ are contained in $X\sm\o K$.
\eenu
\et\vs\noindent

\bt\label{t:7.5}\cite[Theorems 7.7, 7.9]{LiDS}  Let $A\in\sL(X)$ be a positive operator satisfying \eqref{e:2.6}.
Suppose   $K$ is solid, and that   $A$ is weakly irreducible. Then
                \benu
        \item[$(1)$] $r$ is a  simple  eigenvalue of $A$ with  a principal  eigenvector  $w\in \o K$;
                    \item[$(2)$]
                     $\mb{\em GE}_\mu(A)\cap K=\{0\}$ for any  $\mu\in\sig(A)\sm \{r\}$ with $|\mu|>r_e$;
   \item[$(3)$] if $A$ is strongly positive, then
                    \be\label{e:7.1}|\mu|<r,\Hs\A\,\mu\in\sig(A)\sm \{r\}.\ee
           \eenu
\et

\section{Krein-Rutman Type Theorems  for Sectorial Operators}\label{s:3}


Let $X,Y$ be two real Banach spaces with $Y \hookrightarrow X$; moreover, $Y$ is {\em dense} in $X$. Denote by  $\|\.\|$ and $\|\.\|_1$ the norms of $X$ and $Y$, respectively.

Let $K$ be  a cone    in $Y$, 
and  $A$  a closed densely defined operator in $X$ with $-A$ being sectorial (see \cite[Chap. 1]{Henry} for  definition).
We always assume that  the following standing assumptions are fulfilled:
\benu
\item[{\bf (A1)}]\,   $R_\lam(A)K\subset K$ for $\lam>0$ sufficiently large.
\vs
\item[{\bf (A2)}]\, $s_e:=\mb{spb}_e(A)<\mb{spb}(A):=s$.
\item[{\bf (A3)}]\, $\mb{GE}_\mu(A)\subset Y$ for every $\mu\in \sig(A)$ with $\mb{Re}\,\mu> s_e$.
\eenu
As usual, if $s\in\sig(A)$, then we call $s$ the {\em principal eigenvalue} of $A$. Consequently  eigenvectors pertaining to $s$ are referred to as    {\em principal eigenvectors}.

One of our main results is the following  general Krein-Rutman type theorem.

\bt\label{t:8.1} Assume  $K$ is total in $Y$. 
Then the following assertions hold:
 \benu
        \item[$(1)$] $s$ is an  eigenvalue of $A$ admitting   a principal   eigenvector  $u\in K$.

        \item[$(2)$] If \,$\mb{\em int}_YK\ne \emp$ and  contains a principal  eigenvector of $A$, then $s$ shares the same  algebraic and geometric multiplicities.

            \item[$(3)$] All eigenvectors of $A$ corresponding to other  eigenvalues  $\mu\ne s$ with  $\mb{\em Re}\,\mu>s_e$ are contained in $Y\sm\mb{\em int}_YK$.

                             \item[$(4)$] If $\mu\in \sig(A)$, $\mb{\em Re}\,\mu>s_e$ and \,$\mb{\em Im}\,\mu\ne 0$, then
                     \be\label{e:8.1}\mb{\em GE}_\mu(A)\cap K=\{0\}.\ee
           \eenu
\et
{\bf Proof.}  For $t\in\R,\,\,t>s_e$, set  \be\label{e:8.12}\Sigma_1(t):=\{\mu\in\sig(A):\,\,\mb{Re}\,\mu\geq t\}.\ee Since $-A$ is sectorial,     $\Sigma_1(t)$
is a compact subset of $\bbC$ for every  $t>s_e$ (by the definition of sectorial operators). As   every  $\mu\in\sig(A)\sm \sig_e(A)$ is isolated in $\sig(A)$,  one concludes   that $\Sigma_1(t)$ consists of  a finite number of elements.

\vs Let   $\eta\in(s_e,s]$, and  $\Sigma_0(\eta)=\sig(A)\sm \Sigma_1(\eta)$.
By the finiteness of $\Sigma_1(t)$ ($t>s_e$) one trivially checks   that for some  $\de=\de(\eta)>0$,
 \be\label{e:8.9}\ba{ll}\mb{Re}\,\mu\leq\eta-\de,\Hs\A\,\mu\in \Sigma_0(\eta)\ea\ee
Hence  
$\Sigma_0(\eta)$ and $\Sigma_1(\eta)$ form  a spectral decomposition of $\sig(A)$.
Denote by
\be\label{e:8.6}X=X_0(\eta)\oplus X_1(\eta)\ee  the corresponding  decomposition of $X$. Then $X_1(\eta)=\oplus_{\mu\in\Sigma_1(\eta)}\mb{GE}_{\mu}(A)
$ is a  finite-dimensional subspace of $X$. By (A3) we have   $X_1(\eta)\subset Y$.
\vs

For notational simplicity, we rewrite $X_i(\eta):=X_i$ ($i=0,1$).
Let us   split  the  argument below  into several  steps.
\vs

{\bf Step 1.} \,We show that
\be\label{e:8.8}
X_1\cap K\ne\{0\}.
\ee

Let $P=\mb{Cl}_XK$, the closure of $K$ in $X$. Obviously   $P$ is a cone in $X$. Recalling that  $X_1\subset Y$,  to prove \eqref{e:8.8} it suffices to check  that
$$X_1\cap P\ne\{0\}.$$
For this purpose, put $\~A=A-\eta+\de$, where  $\de$ is the positive number in \eqref{e:8.9}.  $\sig(\~A)$ has  a corresponding spectral decomposition $\sig(\~A)=\~\Sigma_0(\eta)\cup\~\Sigma_1(\eta)$ with
$$
\~\Sigma_i(\eta)=\Sigma_i(\eta)-\eta+\de,\Hs i=0,1.$$
We observe that
\be\label{e:8.2}
\sup\{\mb{Re}\,\mu:\,\,\mu\in\~\Sigma_0(\eta)\}\leq-\de,\hs \inf\{\mb{Re}\,\mu:\,\,\mu\in\~\Sigma_1(\eta)\}\geq\de.
\ee
The  direct sum decomposition of $X$ corresponding to the above spectral  decomposition of $\sig(\~A)$ remains  the same as  in \eqref{e:8.6}.

  We claim  that $P\not\subset X_0$. Indeed,  suppose on the contrary  that  $P\subset X_0$. Then since $K$ is total  in $Y$, one would have
$$Y= \mb{Cl}_Y{(K-K)}\subset \mb{Cl}_X{(K-K)}\subset \mb{Cl}_X{(P-P)}\subset X_0.$$ (We emphasize that the closures  $\mb{Cl}_Y$ and $\mb{Cl}_X$  are  taken with respect to the topologies  of $Y$ and $X$, respectively.)
Because $Y$ is dense in $X$, we therefore have $X=\mb{Cl}_XY\subset \mb{Cl}_X X_0=X_0$,  a contradiction.

Take a $u_0\in P\sm X_0$. Write  $u_0=x_0+x_1$, where $x_i\in X_i$. Clearly    $x_1\ne 0$. Let  $u(t)=e^{t\~A}u_0$ ($t\geq 0$), where $e^{t\~A}$ is the $C_0$-semigroup generated by $\~A$.
Then
$$u(t)=e^{t\~A}x_0+e^{t\~A}x_1:=x_0(t)+x_1(t).
$$
We infer from  \eqref{e:8.2}   that
\be\label{e:8.4}
\lim_{t\ra\8}\|x_0(t)\|=0,\hs \lim_{t\ra\8}\|x_1(t)\|=\8.
\ee

By  (A1) we have $R_\lam(A)K\subset K\subset P$ for  $\lam>0$ sufficiently large. Therefore
\be\label{e:8.30}
R_\lam(A)P=R_\lam(A)\,\ol K=\ol{R_\lam(A)K}\subset \ol K\subset P,
\ee
where the closures are taken in $X$. (The second equality in \eqref{e:8.30} is due to the fact that $R_\lam(A)\in\sL(X)$.)
This guarantees that $A$ is {\em semigroup positive}, i.e., $e^{tA}P\subset P$ for $t\geq 0$ (see e.g. Kato \cite[Lemma 5.1]{Kato}). Hence
\be\label{e:8.5}e^{t\~A}P=e^{t(\de-\eta)}e^{tA}P\subset P,\Hs t\geq 0.\ee
In particular, we have  $u(t)=e^{t\~A}u_0\in P$ for all $t\geq 0$.

Now we show that $X_1\cap P\ne\{0\}$ and complete the proof of \eqref{e:8.8}. First, by the first equality  in \eqref{e:8.4} we see that $\lim_{t\ra\8}d\(u(t),X_1\)=0$.
       Now suppose on the contrary that   $X_1\cap P=\{0\}$. Then by \cite[Lemma 2.4]{LiDS} one  deduces that $\lim_{t\ra\8}\|u(t)\|=0$. 
         This contradicts  \eqref{e:8.4}.

\Vs
{\bf Step 2.} \,The verification   of  assertions (1) and (2).
\vs

Take  $\eta=s$. Then
$\Sigma_1(\eta)=\sig_b(A).$ 
 Let $K_1=X_1\cap K$, where $X_i=X_i(\eta)=X_i(s)$ ($i=0,1$) are  given as in \eqref{e:8.6}. Since $X_1$ is a finite-dimensional subspace of $Y$, \eqref{e:8.8} implies  that $K_1$ is a cone in $X_1$.
  As $K$ is total  in $Y$, we have
  $$\ba{ll}
  \mb{Cl}_{X_1}(K_1-K_1)&=\mb{Cl}_{Y}(K_1-K_1)=\mb{Cl}_{Y}(X_1\cap K-X_1\cap K)\\[1ex]
  &\subset \mb{Cl}_{Y}\(X_1\cap (K-K)\)= X_1\cap \mb{Cl}_{Y}(K-K)\\[1ex]
  &=X_1\cap Y=X_1.
  \ea
  $$
  That is,  $K_1$ is total in $X_1$.
 Let  $A_1=A|_{X_1}$. For $\lam>0$ sufficiently large, we infer from \eqref{e:2.6} that $R_\lam(A_1)X_1= X_1$. Thus   by (A1) one easily  verifies  that
 \be\label{e:8.11}R_\lam(A_1)K_1\subset  K_1.\ee

Note that $\sig(A_1)=\Sigma_1(\eta)=\sig_b(A)$.  Let
$$\sig(A_1)=\{\mu_i=s+\bi \b_i:\,\,0\leq i\leq n\}.$$ We may assume that $|\b_0|=\min_{0\leq i\leq n}|\b_i|$.
Fix a number  $\lam>s$ such that \eqref{e:8.11} holds. Then $|\lam-\mu_0|=\min_{0\leq i\leq n}|\lam-\mu_i|$, and hence
  $$r(R_\lam(A_1))=\sup \{1/|\lam-\mu_i|:\,\,0\leq i\leq n\}=1/|\lam-\mu_0|:=r.
  $$
    By Theorem \ref{t:7.4} one concludes  that $r$ is an eigenvalue of $R_\lam(A_1)$ with an eigenvector $w\in K_1$. On the other hand, since $1/|\lam-\mu_i|<r$ for $\mu_i\in\sig(A_1)$ with $\mu_i\ne \mu_0,\ol\mu_0$, we see that  the circle $\bbS_r=\{z\in\bbC:\,\,|z|=r\}$ in the complex plane $\bbC$ contains at most two eigenvalues of  $R_\lam(A_1)$, i.e.,  $1/(\lam-\mu_0)$ and $1/(\lam-\ol\mu_0)$. Thus one necessarily has    $1/(\lam-\mu_0)=1/(\lam-\ol\mu_0)=r$, which implies that $\b_0=0$. It follows that  $\mu_0=s$ is an eigenvalue of $A$; furthermore,   $w$ is an eigenvector of $A$ corresponding to $s$. This  completes the proof of  (1).

If $\mb{int}_YK\ne \emp$ and  contains a principal  eigenvector $v$, one easily verifies that
 $\mb{int}_{X_1}K_1$  is nonvoid and $v\in \mb{int}_{X_1}K_1$.
Thus by Theorem \ref{t:7.4} we deduce  that the  eigenvalue $1/(\lam-s)$ of $R_\lam(A_1)$ has  the same algebraic and geometric multiplicities. Consequently by \eqref{e:2.6} the algebraic and geometric multiplicities of the principal eigenvalue  $s$  of $A$ coincide. Hence assertion (2) holds true.

\Vs
{\bf Step 3.} \,The verification   of  assertions (3) and (4).
\vs

 Let  $\mu\in\sig(A)\sm \{s\}$, $\mb{Re}\,\mu>s_e$.
Take a real number $\eta$ with  $s_e<\eta<s$ such that
$\mu\in \Sigma_1(\eta)$.
Let $X_1=X_1(\eta)$, $K_1=X_1\cap K$, and  $A_1=A|_{X_1}$. Then  as in Step 2 it can be shown that $K_1$ is a total cone in $X_1$. Furthermore, \eqref{e:8.11} remains valid  for $\lam>0$ sufficiently large. Take a  $\lam>s$ such that \eqref{e:8.11} holds and consider the resolvent  operator $R_\lam(A_1)$ of $A_1$ on $X_1$.
Then by Theorem \ref{t:7.4} (4), we deduce  that  $\mb{int}_{X_1}K_1$ does not contain eigenvectors of $R_\lam(A_1)$ pertaining to the  eigenvalue $\lam_\mu:=1/(\lam-\mu)$.  Now if $A$ has an eigenvector $v\in\mb{int}_Y K$ corresponding to $\mu$. Then one easily verify that $v\in \mb{int}_{X_1}K_1$ and is  an eigenvector of $R_\lam(A_1)$ corresponding  to  $\lam_\mu$. This leads to a contradiction and verifies  assertion (3).

\vs If $\mb{Im}\,\mu\ne0$, Theorem \ref{t:7.4} (3) asserts that $\mb{GE}_{\lam_\mu}(R_\lam(A_1))\cap K_1=\{0\}$. We also infer from \eqref{e:2.6} that
$\mb{GE}_{\mu}(A_1)=\mb{GE}_{\lam_\mu}(R_\lam(A_1))$. Therefore  $\mb{GE}_{\mu}(A)\cap K=\mb{GE}_{\mu}(A_1)\cap K_1=\{0\}$. This  completes the proof of assertion (4). $\bx$

\br Note that we do   not require that $e^{tA}Y\subset Y$ for $t\geq 0$ in the proof of the above theorem. This allows us to avoid deriving higher regularity results on the corresponding parabolic equations when applying the theory to  elliptic differential operators.

\er
\bt\label{t:8.2} In addition to  $(A1)-(A3)$, we also  assume that the following  hypothesis is fulfilled:
\benu\item[{\em{\bf(A4)}}] $K$ is a solid cone in $Y$, and
$R_\lam(A)(K\sm\{0\})\subset \mb{\em int}_YK$.
\eenu Then  $s$ is a simple  eigenvalue with  a unique normalized  eigenvector  $w\in \mb{\em int}_YK$. Moreover, for any  $\mu\in\sig(A)\sm \{s\}$ with $\mb{\em Re}\,\mu>s_e$,
                     \be\label{e:8.15}\mb{\em GE}_\mu(A)\cap K=\{0\}.\ee
\et
{\bf Proof.} The proof follows a fully analogous manner as the one for Theorem \ref{t:7.5}, and is thus omitted. The interested reader may consult \cite[Theorems 7.3, 7.7]{LiDS} for details. $\bx$

\br
Under the hypotheses of Theorem \ref{t:8.2}, one may expect that the boundary spectrum $\sig_b(A)$ consists of  exactly one eigenvalue of $A$. Unfortunately the easy counterexample below  indicates that this  may  fail to be true.
\er\noindent
{\bf Example 3.1.} Let $X=Y=\R^3$. For computational convenience, here we make a convention that  $\R^3$ consists of column vectors. Denote $v'$  the transpose of a row vector $v=(x,y,z)$.
Define a cone $K$  in $X$ as
$$
K=\{(x,y,z)'\in X:\,\,z\geq \sqrt{x^2+y^2}\}.
$$
Then $\mb{int}_YK=\{(x,y,z)'\in X:\,\,z> \sqrt{x^2+y^2}\}$. Let
$$
A=\(\begin{matrix}0&-1&0\\1&0&0\\0&0&0\end{matrix}\).
$$
For $\lam>0$, simple computations yield
$$
R_\lam(A):=(\lam-A)^{-1}=\(\begin{matrix}B& O\\[1ex]
O&{\lam^{-1}}\end{matrix}\),\hs \mb{where }\,B=\frac{1}{1+\lam^2}\(\begin{matrix}\lam&-1\\1&\lam\end{matrix}\).
$$
 Therefore
$$
R_\lam(A)u=\frac{1}{1+\lam^2}\(\lam x-y,\,x+\lam y,\, \frac{1+\lam^2}{\lam}z\)':=\frac{1}{1+\lam^2}(\~x,\~y,\~z)'.
$$
Observe that
\be\label{e:8.16}
\~x^2+\~y^2=(x^2+y^2)+\lam^2(x^2+y^2)
=(1+\lam^2)(x^2+y^2).
\ee

Now let $u=(x,y,z)'\in K$. Then $x^2+y^2\leq z^2$. Since $(1+\lam^2)/\lam^2>1$, by \eqref{e:8.16} we deduce that
$$\~x^2+\~y^2\leq (1+\lam^2)z^2\leq \(\frac{1+\lam^2}{\lam}\)^2z^2=\~z^2.$$
This implies $(\~x,\~y,\~z)'\in K$. Thus we see that $R_\lam(A)K\subset K$.

If $u\in \pa K$, $u\ne 0$, then $x^2+y^2=z^2\ne0$. By \eqref{e:8.16} we find that
$$
\~x^2+\~y^2=(1+\lam^2)z^2<\(\frac{1+\lam^2}{\lam}\)^2z^2=\~z^2.
$$
Hence $(\~x,\~y,\~z)'\in \mb{int}_YK$.  Therefore  $R_\lam(A)u\in \mb{int}_YK$. This indicates that the operator given by $A$ satisfies all the requirements in Theorem \ref{t:8.2}. However, all the eigenvalues  of $A$ has the same real part $s=0$.
\Vs
To guarantee the uniqueness of elements in $\sig_b(A)$,  Nussbaum \cite{Nuss2} used the notion of ``$u_0$-positivity\,'' due to Krasnosel'skii \cite{Kras}; see \cite[Theorem 1.3]{Nuss2}. Here we remark  that if the  semigroup $e^{tA}$  has some strong positivity property, then one can still ensure the uniqueness of elements in $\sig_b(A)$.

\bt\label{t:8.3}In addition to  $(A1)-(A4)$, assume  that
\benu\item[{\em{\bf(A5)}}] for any $t>0$ and $\mu\in\sig(A)$ with  $\mb{\em Re}\,\mu>s_e$,
$$
e^{tA}(K_\mu\sm\{0\})\in \mb{\em int}_YK,\hs \mb{where }K_\mu=\mb{\em GE}_\mu(A)\cap K.$$
\eenu
 Then $\sig_b(A)=\{s\}$. 
\et
{\bf Proof.} Let $Y'=\oplus_{\mu\in\sig_b(A)}\mb{GE}_\mu(A)$, $A'=A|_{Y'}$. Denote $K'=Y'\cap K$. Then by Theorem \ref{t:8.2} we see that $K'\ne\{0\}$. Hence $K'$ is a cone in $Y'$. Let $v\in K'$, $v\ne 0$. By (A5) we have  $e^{tA}v\in \mb{int}_YK$ for $t>0.$ Since $e^{tA}v\in Y'$,  one  trivially verifies   that $e^{tA}v\in \mb{int}_{Y'} {K'}$. Therefore  $K'$ is a solid cone in $Y'$ and \be\label{e:8.13}e^{tA'}(K'\sm\{0\})\subset \mb{int}_{Y'} {K '},\Hs t>0.\ee

Now let $\mu:=s+\bi\b\in\sig_b(A)$. Then $\lam :=e^{\mu  t}=e^{st}e^{\bi\b t}$ is an eigenvalue of $e^{tA'}$ with $|\lam |=e^{st}:=r(t)$. But \eqref{e:8.13} implies that $r(t)$ is the unique eigenvalue of $e^{tA'}$ on the circle $\bbS_{r(t)}:=\{z\in\bbC|\,\,|z|=r(t)\}$ for $t>0$. Hence we necessarily have  $e^{\bi\b t}=1$, and therefore $\b t\in \{2k\pi:\,\,k\in \mathbb{Z}\}$ for all  $t>0$. But this is impossible unless $\b =0$. This proves what we desired. $\bx$

\section{Principal Eigenvalue Problem of Elliptic  Operators Associated with Degenerate Mixed Boundary Conditions}\label{s:4}
\vs
As an illustrating example,  we consider  the principal eigenvalue problem of   the elliptic  operator $L$ on  a smooth bounded domain $\Omega\subset\mathbb{R}^n(n\geq 1)$:  
\begin{equation}\label{equ:1-2}
L u=-\sum^n_{i,j=1}a_{ij}(x)\frac{\partial^2u}{\partial x_i\partial x_j}+\sum^n_{i=1}b_i(x)\frac{\partial u}{\partial x_i}+c(x)u,
\end{equation}
which is  associated with the  mixed boundary condition:
\be\label{e:bc}Bu:=\a(x)u+\b(x)\frac{\d u}{\d \nu}=0\ee
on the boundary $\Gamma:=\d\W$ of $\W$, where   $\nu$ stands for  the unit outward normal vector field on $\Gamma$. The coefficients of $L$ and $B$ 
are assumed to be $C^\8$ functions satisfying  the hypotheses below:
\benu
\item[{\bf(H1)}]
$a_{ij}=a_{ji}$ ($1\leq i,j\leq n$), and   there is $\theta>0$ such that
$$
\sum^n_{i,j=1} a_{ij}(x)\varsigma_{i}\varsigma_j \geq\theta|\varsigma|^2, \Hs \A\,\varsigma\in\R^n,\,\,x\in \bar{\W};
$$
\item[{\bf(H2)}]  $c,\a,\b$ are nonnegative functions satisfying that
\be\label{e:8.14}\a(x)+\b(x)>0,\Hs\A\,x\in\Gamma.\ee
\eenu
In the case of the Dirichet boundary condition or the  Robin boundary condition (regular case), this problem has already been well understood; see e.g. \cite[Theorem 12.1]{Amann}, \cite[Theorem 1.4]{Du} and also \cite{GD, Nuss2}. Here we are interested in a degenerate case where $\b$ may vanish on a proper subset of $\Gamma$. In such a situation, if  $L$ has a divergence form (hence $L$ enjoys some symmetric properties), one can find some nice results concerning the principal eigenvalue problem  of $L$  in  Taira \cite[Theorem 1.2]{Taira}. As an  application of the abstract results given in Section \ref{s:3},   we deal with  the general case and present a less involved argument on the  problem.

\vs
\noindent{\bf $\bullet$ Some fundamental results.} \,First, making use of  the classical Hopf's lemma, one can easily verify the comparison result  below:
\bl\label{l:8.2}
Let $u\in C^1(\ol\W)\cap C^2(\W)$, $u\not\equiv 0$. Assume that
$$Lu+\lam u\geq0 \,\,(\mb{in }\W),\hs Bu\geq0\,\,(\mb{on }\Gamma),$$
where $\lam\geq 0$. Then $u(x)>0$ for $x\in\W$.
\el

Denote by $W^{s,p}(\W)$ ($s\in \R_+$, $1\leq p<\8$) the Sobolev spaces equipped with the standard norms. We infer from
Taira \cite[pp.5, Theorem 1]{Taira} that the following existence and uniqueness result holds true.

\bl\label{l:8.1} Let $1<p<\8$, $s>1+1/p$, and let $\lam\geq 0$. Then for any $f\in W^{s-2,p}(\W)$, the homogeneous boundary value problem
\be\label{e:8.17}Lu+\lam u=f\mb{ $($in\,\,$\W$$)$},\hs \mb{ $Bu=0$ $($on\,\,$\Gamma$$)$}\ee has a unique solution $u\in W^{s,p}(\W)$. Here the boundary condition  is understood in the sense that $B$ can be viewed as a linear operator  from $W^{s,p}(\W)$ to Besov space $B_*^{s-1-1/p,\,p}(\Gamma)$ (see \cite{Taira} pp. 3 for details).
\el

Note that  Lemma \ref{l:8.1} implies   that if $f\in C^1(\ol\W)$, then the solution $u$ of \eqref{e:8.17} belongs to $C^2(\ol\W)$, and hence it solves \eqref{e:8.17} in the classical sense. Indeed, if $f\in C^1(\ol\W)$ then $f\in W^{1,p}(\W)$ for any $1<p<\8$. Lemma \ref{l:8.1} then asserts that $u\in W^{3,p}(\W)$. Taking  a number  $p>1$ sufficiently large so that  $W^{3,p}(\W)\hookrightarrow C^{2}(\ol\W)$, one immediately concludes that $u\in C^{2}(\ol\W)$.

 By virtue of \cite[pp.4, Theorem 1]{Taira} we also  deduce   that
\be\label{e:8.18b}\|u\|_{C^{2}(\ol\W)}\leq C \|u\|_{W^{3,p}(\W)}\leq C \|f\|_{W^{1,p}(\W)}\leq C\|f\|_{C^{1}(\ol\W)}\ee
for all $f\in C^1(\ol\W)$, where $C$ denotes a general constant independent of $f$.
\vs

\noindent{\bf $\bullet$ Resolvent strong positivity of  $A=-L$.} \,Let $X=L^2(\W)$, and set
$$Y=\{u\in C^1(\ol\W):\,\,\mb{$u$ satisfies \eqref{e:bc}}\}.$$  $Y$ is equipped  with the usual norm of $C^1(\ol\W)$. Clearly $Y\hookrightarrow X$.

Let $K$ be the positive cone in $Y$ consisting of nonnegative functions.

Denote by $A$ the operator $-L$ with domain
$$\cD(A)=\{u\in H^2(\W):\,\,Bu=0\},$$ where the boundary condition $Bu=0$ is understood in the same sense as in Lemma \ref{l:8.1}. Invoking \cite[pp.5, Theorem 2]{Taira} we deduce  that   $-A$ is a sectorial operator in $X$ with compact resolvent. 
Hence by \cite[Ex. 2.4, (i)]{Kry} it is easy to see that
$
\mb{spb}_e(A)=-\8<\mb{spb}(A)<\8.
$
Thus $A$ fulfills  (A2) in Section \ref{s:3}.

 We infer from Lemma \ref{l:8.1} and \eqref{e:8.18b} that $R_\lam(A)Y\subset Y$ for $\lam\geq 0$; furthermore,   $R_\lam(A)|_{Y}$ is compact as an operator on the space $Y$.

The following  result  indicates  that $A$ fulfills  hypotheses (A1) and (A4) in Section \ref{s:3}. The proof of such a result is   somewhat  standard. We include the details  in the Appendix part  for the readers' convenience.
 \bl\label{l:8.3}
$R_\lam(A) (K\sm\{0\})\subset \mb{\em int}_YK$ for each $\lam\geq 0$.
\el

\noindent{\bf $\bullet$ Regularity of the generalized eigenfunctions.}\,
Let $\bbA$ be the complexification of $A$ with  $\cD(\bbA)=\cD(A)+\bi\cD(A)$.
 We start with  the  eigenfunctions of $\bbA$. Let $\mu=a+\bi b\in\sig(A)$, and let $w=u+\bi v$ be a corresponding eigenfunction of $\bbA$, where $u,v\in D(A)$. Then
$\bbA w=\mu w$ amounts to say that
\be\label{e:8.25}
Au=au-bv,\hs Av=av+bu.
\ee
Since $u,v\in  H^2(\W)$, by Lemma \ref{l:8.1} and \eqref{e:8.25} one finds  that $u,v\in H^4(\W)$. Further by a simple bootstrap argument we finally conclude that $u,v\in H^s(\W)$ for all $s\geq 0$. It follows by the Sobolev embeddings that $u,v\in C^\8(\ol\W)$.

Now let $g\in\mb{GE}_\mu(\bbA)$ and  $\mb{rank}\,(g)\geq 2$. Set $k=\mb{rank}\,(g)-1$. Then $(\bbA-\mu)^kg:=w$  is an eigenfunction of $\bbA$. Hence $w$ a $C^\8$ function on $\ol\W$.

Note that it is readily implied in $(\bbA-\mu)^kg =w$ that $(\bbA-\mu)^jg\in \cD(\bbA)$ for all $j\leq k$. In particular,
$$(\bbA-\mu)^{k-1}g:=f_1\in\cD(\bbA)\subset \bbH^2(\W):=H^s(\W)+\bi H^s(\W).$$
Therefore by $(\bbA-\mu)f_1=w$ we find that
\be\label{e:8.28}
\bbA f_1=w+\mu f_1:=\~f_1\in \bbH^2(\W).
\ee
It follows by Lemma \ref{l:8.1} that $f_1\in \bbH^4(\W)$. This in turn  implies that $\~f_1\in \bbH^4(\W)$. \eqref{e:8.28} and  Lemma \ref{l:8.1} then asserts that $f_1\in \bbH^6(\W)$. Once again by a bootstrap argument we see  that  $f_1\in \bbH^s(\W)$ for all $s\geq 0$.

Repeating the above argument with $w$ and $f_1$ therein  replaced by $f_1$ and $(\bbA-\mu)^{k-2}g:=f_2$, respectively, one deduces that $f_2\in \bbH^s(\W)$ for $s\geq 0$. Continuing this procedure we finally obtain that $g=f_k\in \bbH^s(\W)$ for all $s\geq 0$. The Sobolev embeddings then immediately imply   that $g$ is a  $C^\8$ function.

It follows from the above results that $\mb{GE}_\mu(A)\subset Y$ for any $\mu\in\sig(A)$. Hence  $A$ fulfills hypothesis (A3).
\vs
By far we have seen that the operator   $A$ satisfies hypotheses (A1)-(A4).
\vs

\noindent{\bf $\bullet$ The verification of  hypothesis (A5).}\,
Let $\mu\in\sig(A)$. Denote by $\bbA_\mu$ the restriction of $\bbA$ on $\mb{GE}_\mu(\bbA)$. Given  $g\in \mb{GE}_\mu(\bbA)$, 
let $u=u(t)$ be the solution of  equation $\dot u=\bbA u$ with  $u(0)=g$. Then $u(t)=e^{t\bbA}g=e^{t\bbA_\mu}g$. Since  $\bbA_\mu$ is a bounded operator on $\mb{GE}_\mu(\bbA)$,   we have
\be\label{e:8.26}
\ba{ll}u(t)&=e^{t\bbA_\mu}g=e^{\mu t}e^{t(\bbA_\mu-\mu)}g=e^{\mu t}\sum_{j=0}^\8\frac{t^j}{j!}(\bbA_\mu-\mu)^jg\\[1ex]
&=e^{\mu t}\(I+\frac{t}{1!}( \bbA_\mu-\mu)+\cdots+\frac{t^{k-1}}{(k-1)!}(\bbA_\mu-\mu)^{(k-1)}\)g,
\ea
\ee
where $k=\mb{rank}\,(g)$. Noticing that $(\bbA_\mu-\mu)^{j}g\in \mb{GE}_\mu(\bbA)$ for any integer $j\geq 0$, by what we have proved above it is clear   that $(\bbA_\mu-\mu)^{j}g$ is a $C^\8$ function. Consequently for each $t\geq 0$ fixed, $u(t)$ is a $C^\8$ function in the  space variable on $\ol\W$.
We write  $u(t,x)=u(t)(x)$ for $(t,x)\in \R_+\X\ol\W$. Then it can be easily seen that $u=u(t,x)$ is a complex   $C^\8$ function on $\R_+\X\ol\W$.

Now we come back to the real situation. The above result  implies  that for each  $v_0\in \mb{GE}_\mu(A)$, the function  $v(t,x):=v(t)(x)$, where $v(t)=e^{tA}v_0$,  is a $C^\8$ function on $\R_+\X\ol\W$.
Therefore  $v$ is a  classical solution of the   parabolic equation:
\be\label{e:8.27}
\frac{\pa v}{\pa t}+L v=0,\Hs x\in\W,\,\,t>0
\ee
associated with the boundary condition $Bv=0$ on $\Gamma$. Thanks to the Hopf's lemma for parabolic equations (see e.g. Friedman \cite[Chap. 2, Theorem 14]{Friedman} or Smith \cite[Chap. 7, Theorem 2.2]{Smith}), using almost the same argument as in the proof of Lemma \ref{l:8.3} (see the Appendix), it can be shown that $v(t,\.)\in \mb{int}_YK$ if $v_0\in K$, $v_0\ne0$. This is precisely what we desired.

\Vs
Now   that $A$ satisfies hypotheses $(A1)-(A5)$, as a straightforward application of Theorems \ref{t:8.2} and \ref{t:8.3}, one immediately obtains the following result.
\bt\label{t:8.5} The following assertions hold true:
\benu\item[$(1)$] The spectral bound $s$ of $A$ is an algebraically simple  eigenvalue  with  a corresponding  eigenvector  $w\in \mb{\em int}_YK$.
\item[$(2)$] $\mb{\em GE}_\mu(A)\cap K=\{0\}$ for any  $\mu\in\sig(A)\sm \{s\}$.
\item[$(3)$] $\sig_b(A)=\{s\}$
   \eenu
\et

\br
We mention that in the case of the Dirichlet 
(or  Robin) boundary condition, 
almost all the facts concerning the operator $A$ needed in proving  Theorem \ref{t:8.3} are well known and need not be checked. Therefore the theorem  becomes  nearly  an immediate consequence of Theorems \ref{t:8.2} and \ref{t:8.3}.\er

\Vs\Vs
\centerline{\bf Appendix: Proof of Lemma \ref{l:8.3}}
\Vs
\noindent{\bf Proof.}
Let $f\in K\sm\{0\}$, and $u=R_\lam(A) f$. Then  $u\in C^2(\ol\W)$. By Lemma \ref{l:8.2} we deduce that $u(x)>0$ for $x\in \W$. Set $\Gamma_0=\{x\in\Gamma:\,\,u(x)=0\}$. Then  the classical  Hopf's lemma asserts  that $\frac{\pa u}{\pa \nu}(x)<0$ for $x\in\Gamma_0$. Hence by compactness of $\Gamma_0$,  there is $\ve_0>0$ such that
$\frac{\pa u}{\pa\nu}(x)\leq -3\ve_0$ for $x\in \Gamma_0$.

Denote $\|\.\|_1=\|\.\|_{C^1(\ol\W)}$. Take  a neighborhood $\cW$ of $\Gamma_0$ in $\Gamma$ such that
\be\label{e:8.18}
\frac{\pa u}{\pa\nu}(x)\leq -2\ve_0,\Hs x\in \cW.
\ee
Then there exists    $\de>0$ such that for all $h\in Y$ with $\|h\|_1<\de$,
\be\label{e:8.19}
\frac{\pa (u+h)}{\pa\nu}(x)\leq -\ve_0<0,\Hs x\in \cW.
\ee

We claim that $\a(x)>0$ for $x\in\Gamma_0$. Indeed, if $\a(x)=0$ then by \eqref{e:8.14} we have $\b(x)>0$. Thus by  \eqref{e:8.18} one deduces that $Bu(x)=\b(x)\frac{\pa u}{\pa\nu}(x)\ne 0$, a contradiction. Hence the claim holds true. By compactness of $\Gamma_0$ we deduce that $\a(x)\geq 2\ve_1>0$ for all $x\in\Gamma_0$. Therefore by continuity of $\a$ it can be assumed  that the neighborhood $\cW$ of $\Gamma_0$ is chosen sufficiently small so that
\be\label{e:8.20}
\a(x)\geq\ve_1>0,\Hs x\in\cW.
\ee

Now  for any $h\in Y$ with $\|h\|_1<\de$, we have at any point $x\in\cW$ that
$$
\a(x)(u+h)(x)=-\b(x)\frac{\pa (u+h)}{\pa\nu}(x)\geq(\mb{by }\eqref{e:8.19})\geq 0.
$$
Hence by \eqref{e:8.20} we see that
\be\label{e:8.21}
(u+h)(x)\geq 0,\Hs x\in \cW.
\ee
Using \eqref{e:8.19} and \eqref{e:8.21} it is not difficult  to deduce that
there is a neighborhood $U$ of $\Gamma_0$ in $\ol\W$ such that for any $h\in Y$ with $\|h\|_1<\de$,
\be\label{e:8.22}
(u+h)(x)\geq 0,\Hs x\in U.
\ee

We may assume that $U$ is open relative to $\ol\W$. Hence $\Gamma_1:=\Gamma\sm U$ is compact. Because $u$ is positive on $\Gamma_1$, there is $\ve_2>0$ such that $u(x)\geq 2\ve_2$ for $x\in \Gamma_1$. This allows us to pick a neighborhood $V$ of $\Gamma_1$ in $\ol\W$ such that
$u(x)\geq \ve_2$ for $x\in V$. Further one can restrict $\de$ sufficiently small so that
\be\label{e:8.23}
(u+h)(x)\geq 0,\Hs x\in V
\ee
for all $h\in Y$, $\|h\|_1<\de$. Note that $G:=U\cup V$ is a neighborhood of $\Gamma$ in $\ol\W$.

It can be assumed that both $U,V$ are open relative to $\ol\W$. Hence $G$ is open in $\ol\W$. Consequently $F:=\ol\W\sm G$ is a compact subset of $\ol\W$.
Since $u$ is positive on $F$, there exists $\ve_3>0$ such that $u(x)\geq 2\ve_3$ for $x\in F$. Therefore if $\de$ is sufficiently small then $u+h$ is positive on $F$  for all $h\in Y$ with $\|h\|_1<\de$. Combining this with  \eqref{e:8.22} and \eqref{e:8.23} it follows  that
$u+h\geq 0$ in $\ol\W$ for all $h\in Y$ with $\|h\|_1<\de$, i.e., $u+h\in K$. Hence $u\in\mb{int}_YK$. $\bx$
\Vs



\end{document}